\documentclass[12pt,reqno]{amsart}
\usepackage{amssymb}  
\usepackage{pifont}  
\usepackage{marvosym}  

\input cyracc.def
\newfam\cyrfam

\font\twelvecyr=wncyr10 scaled 1200 
\def\cyr{\fam\cyrfam\twelvecyr\cyracc}

\renewcommand\AA{{\mathcal A}}
\newcommand\BB{{\mathcal B}}
\newcommand\MM{{\mathcal M}}

\newcommand\RRR{{\mathbb R}}
\newcommand\NNN{{\mathbb N}}
\newcommand\DDD{{\mathbb D}}

\newcommand\TTT{{\mathbb T}}
\newcommand\CCC{{\mathbb C}}

\newcommand\cchi{{\raise 2 pt \hbox{$\chi$}}}

\newcommand\rest{\upharpoonright}     
\newcommand\res{\mathord {\upharpoonright}}  

\newcommand\JMP{\mathsf{JMP}} 

\newcommand\STEP{\mathsf{STEP}} 
\newcommand\SF{\mathsf{SF}} 

\renewcommand\Re{\mathrm{Re}}   
\newcommand\st{\mathrm{st}}   

\newcommand\supt{\mathrm{supt}}   

\newcommand\Sh{\mbox{\cyr Sh}}  

\newcommand\ltn{\mathord{\!\! \lessdot}}  

\newcommand\lea{\sqsubseteq}  

\newcommand\iv{^{-1}} 

\newenvironment{itemizz}{\begin{itemize}}  {\end{itemize}}                              

\newenvironment{itemizn}[1] 
{\begin{itemize} %
\renewcommand\labelitemi{\ding{#1}}} %
{\end{itemize}}

\newtheorem{theorem}{Theorem}[section]
\newtheorem{definition}[theorem]{Definition}
\newtheorem{lemma}[theorem]{Lemma}
\newtheorem{corollary}[theorem]{Corollary}
\newtheorem{proposition}[theorem]{Proposition}
\newtheorem{example}[theorem]{Example}

\begin{document}

\title{The Complex Stone--Weierstrass Property}
\author{Kenneth Kunen}
\thanks{Author partially supported by NSF Grant DMS-0097881.}
\address{Department of Mathematics\\
University of Wisconsin\\
Madison, WI 57306, USA}
\email{kunen@math.wisc.edu}
\urladdr{http://www.math.wisc.edu/\symbol{126}kunen}
\date{\today}

\begin{abstract}
The compact Hausdorff space $X$ has the CSWP iff
every subalgebra of $C(X, \mathbb{C})$ which 
separates points and contains the constant functions
is dense in $C(X, \mathbb{C})$. 
Results of W. Rudin (1956), and Hoffman and Singer (1960), show that
all scattered $X$ have the CSWP and many non-scattered $X$ fail
the CSWP, but it was left open whether having the CSWP is just equivalent
to being scattered.

Here, we prove some general facts
about the CSWP; and in particular we show that if $X$ is
a compact ordered space, then $X$ has the CSWP
iff $X$ does not contain a copy of the Cantor set.
This provides a class of non-scattered spaces with the CSWP.
\end{abstract}
\maketitle

\section{Introduction} 
\label{sec-intro}
All spaces discussed in this paper are Hausdorff.

\begin{definition}
If $X$ is compact, then
$C(X) = C(X,\CCC)$ is the algebra of continuous
complex-valued functions on $X$,
with the usual supremum norm.
$\AA \lea C(X)$ means that $\AA$ is a subalgebra of $C(X)$
which separates points and contains the constant functions.
\end{definition}

If $\AA \lea C(X)$ is self-adjoint
($f \in \AA \leftrightarrow \overline f \in \AA$),
then $\AA$ is dense in $C(X)$ by the standard
Stone-Weierstrass Theorem for real-valued functions.

\begin{definition}
The compact space $X$
has the \emph{Complex Stone--Weier\-strass Property (CSWP)}
iff every $\AA \lea C(X)$ is dense in $C(X)$.
\end{definition}

The CSWP may be considered a notion of ``smallness'' by the following easy

\begin{lemma}
\label{lemma-cswp-down}
If $X$ is a closed subspace of the compact space $Y$,
and $Y$ has the CSWP, then $X$ has the CSWP.
\end{lemma}
\begin{proof}
If $\AA \lea C(X)$, let $\AA' = \{f \in C(Y) : f\res X \in \AA\}$.
If $\AA'$ is dense in $C(Y)$,
then $\AA$ is dense in $C(X)$.
\end{proof}

Classical examples show that the CSWP can fail.  As usual,
$D = \{z\in\CCC : |z| < 1\}$ and 
$\TTT = \{z\in\CCC : |z| = 1\}$.  Let $\AA \lea C(\overline D)$ be
the \textit{disc algebra}; these are the functions in $C(\overline D)$
which are holomorphic on $D$.
Then $\AA$ shows that $\overline D$ fails to have the CSWP,
and $\AA\rest\TTT$ shows that $\TTT$ fails to have the CSWP\@.
Here, the restricted algebra $\AA\rest\TTT$  is defined by:

\begin{definition}
Given $\AA \lea C(X)$ and a closed subset $H \subseteq X$,
let $\AA\rest H = \{f \res H : f \in \AA\} \lea C(H)$.
\end{definition}

W. Rudin proved the following basic results about the CSWP:

\begin{theorem}[Rudin]
\label{thm-rudin} 
For all compact $X$:
\begin{itemizz}
\item[1.] If $X$ contains a copy of the Cantor set,
then $X$ does not have the CSWP.
\item[2.] If $X$ is scattered, then $X$ has the CSWP.
\end{itemizz}
\end{theorem}

As usual, $X$ is scattered iff it has no perfect subsets;
see also Definition \ref{def-kernel} below.
(1) is proved in \cite{RUD1} and 
(2) is proved in \cite{RUD2}.
Of course, for (1), it is sufficient (by Lemma \ref{lemma-cswp-down})
to prove that the Cantor set itself fails the CSWP, and this is done
by another use of algebras of holomorphic functions.

This theorem completely characterizes the CSWP for compact metric $X$,
since for such spaces, $X$ is scattered iff $X$ does not
contain a copy of the Cantor set.  
However, there are non-metrizable non-scattered compact spaces,
such as $\beta\NNN$, which do not contain Cantor subsets.

We do not know a simple characterization of the CSWP which
holds for all compact $X$, but by the results in this paper
and earlier results of Hoffman and Singer \cite{HS},
neither of the implications (1) and (2) reverses.
(1) is not an ``iff'', since
$\beta \NNN$ does not have the CSWP by \cite{HS};
see Section \ref{sec-ex} for more on their example.
(2) also is not an ``iff'' by the following
result, which shows that (1) \textit{is} an ``iff'' for a restricted
class of spaces:

\begin{theorem}
\label{thm-lots}
Let $X$ be totally ordered by $<$, and assume that $X$
is compact in its order topology.  Then $X$ has the CSWP
iff $X$ does not contain a copy of the Cantor set.
\end{theorem}

In particular, Aleksandrov's double arrow space, which is not scattered,
has the CSWP,
since it does not contain a copy of the Cantor set.
The double arrow space is usually obtained by replacing
each point in $[0,1]$ by a pair of neighboring points (see
Section \ref{sec-lots}).
It may also be viewed as the maximal ideal space of
the closed algebra generated by the piecewise continuous functions
on $\TTT$ (see Lemma \ref{lemma-double-arrow}).
With this identification, the CSWP for the double arrow space
may be interpreted as a statement about algebras of piecewise continuous
functions (see Corollary \ref{cor-jump}).

Of course, one direction of Theorem \ref{thm-lots} is contained in
Theorem \ref{thm-rudin}.
We prove the other direction in the case that $X$ is
\textit{separable} in Section \ref{sec-lots},
which also lists some elementary properties of compact ordered spaces.
In particular (see Lemma \ref{lemma-disc}),
if such an $X$ is separable and fails to contain a Cantor set,
then it must be totally disconnected;
equivalently (since $X$ is compact), the subsets of $X$ which
are \textit{clopen} (both closed and open) form a base for the topology.
If $H$ is clopen, then its characteristic function $\cchi_H \in C(X)$
is an idempotent.  Section \ref{sec-id} contains further results
on idempotents; these results reduce the proof
of the separable case of Theorem \ref{thm-lots} to showing that for such $X$,
every closed $\AA \lea C(X)$ contains a non-trivial idempotent.

The proof of Theorem \ref{thm-lots}, without the restriction that
$X$ be separable, is given in Section \ref{sec-meas}.
This is actually quite easy, using the observation
(Lemma \ref{lemma-supt-cswp})
that if $X$ fails the CSWP, then there is a $Y \subseteq X$
which also fails the CSWP such that $Y$ is the (closed)
support of a measure.  In the case of compact ordered
spaces, it is already known that the support of a measure is separable.

We remark that it is easy to construct a connected \textit{non}-separable
compact ordered space $X$ which does not contain a Cantor set.
It is curious that the proof of the
CSWP for $X$ proceeds via the study of idempotents, even
though $C(X)$ itself has only the trivial
idempotents.

Sections \ref{sec-ker} considers $\AA \rest \ker(X)$, where
$\ker(X)$ is defined by:

\begin{definition}
\label{def-kernel}
For a topological spaces $X$:
\begin{itemizz}
\item[1.]
$Y \subseteq X$ is \emph{perfect} iff $Y$ is non-empty and closed
and no point of $Y$ is isolated in $Y$.
\item[2.]
$X$ is \emph{scattered} iff $X$ has no perfect subsets.
\item[3.]
The \emph{kernel}, $\ker(X)$, is $\emptyset$ if $X$ is scattered,
and the largest perfect subset of $X$ otherwise.
\end{itemizz}
\end{definition}

In (3), one constructs the largest perfect subset
by taking the closure of the union of all perfect subsets.
For a compact non-scattered $X$,
all idempotents
of $\AA \res \ker(X)$ extend to idempotents of $\AA$
(see Lemma \ref{lemma-restrict-kernel}), and $X$
 has the CSWP iff $\ker(X)$ has the CSWP
(see Corollary \ref{cor-CSWP-iff}).

Some further remarks on these restrictions $\AA\res H$:
Our definition of the property ``$\AA \lea C(X)$''
did not require that $\AA$ be closed in $C(X)$.
Of course, to verify the CSWP for a given $X$,
one need only consider closed $\AA$.  
However, we shall frequently be studying restrictions of $\AA$.
Even if $\AA$ is closed, $\AA\res H $ need not be
closed in $C(H)$.
In fact, by Glicksberg \cite{GL}, if $\AA \lea C(X)$ 
and $\AA\res H$ is closed in $C(H)$ for all closed $H \subseteq X$,
then $\AA = C(X)$.
$\AA\res H$ \textit{is} closed in $C(H)$ if
$H$ contains the \v Silov boundary (see Section \ref{sec-id}),
or if $H = \ker(X)$ (see Lemma \ref{lemma-restrict-kernel}).

If $\AA$ is a commutative Banach algebra with a unit element,
we shall use $\MM(\AA)$ to denote the
\textit{maximal ideal space} or \textit{spectrum}.
Elements of $\MM(\AA)$
may be viewed either as maximal ideals in $\AA$, or as homomorphisms
from $\AA$ to $\CCC$. 
The Gel'fand transform maps $\AA$ into a subalgebra of
$C(\MM(\AA))$.  For more on these notions, 
see, e.g., \cite{GA,HO,RUD4}.

One may also define a ``$\lea$'' relation between 
commutative unital Banach algebras: $\AA \lea \BB$ iff
$\AA$ is a subalgebra of $\BB$ and $\AA$ separates
the points of
$\MM(\BB)$ (mapping $\AA, \BB$ into subalgebras of $C(\MM(\BB))$).
One might say that $\BB$ has the CSWP iff $\AA \lea \BB$
implies $\AA = \BB$.
This is related to the Stone-Weierstrass Property (SWP) discussed
by Katznelson and Rudin \cite{KR}; $\BB$ has the SWP iff every
\textit{self-adjoint} $\AA \lea \BB$ equals $\BB$.

\section{Idempotents and Restrictions}
\label{sec-id}
We begin with a useful criterion for telling whether $\AA \res H$
is closed in $C(H)$.

\begin{definition}
$\|f\|_H = \sup\{|f(x)| : x \in H\}$.
\end{definition}

\begin{lemma}
\label{lemma-closed-equiv}
Assume that $\AA\lea C(X)$ and is closed in $C(X)$.
Let $H$ be a closed subset of $X$.  Then the following are
equivalent:
\begin{itemizz}
\item[1.]
$\AA\res H$ is closed in $C(H)$.
\item[2.]
There is a finite constant $c \ge 1$ such that for all $f \in \AA$,
there is an $f^* \in \AA$ with $f \res H =  f^*  \res H$ and
$\| f^*\| \le c \|f\|_H$.
\end{itemizz}
\end{lemma}
\begin{proof}
As pointed out in \cite{GL}, $(1) \to (2)$ follows by applying
the Open Mapping Theorem to the restriction map from $\AA$ onto $\AA \res H$.

For $(2) \to (1)$, we repeat one of the steps in the proof 
of the Open Mapping Theorem.
Suppose that we have $f_n \in \AA$ for $n \in \NNN$,
with $f_n \res H$ converging in $C(H)$.  We must find
a $g \in \AA$ such that $g\res H = \lim_n (f_n \res H)$.
WLOG, we may assume that each $\|f_{n+1} - f_{n}\|_H \le 2^{-n}$.
For each $n$, apply (2) and choose $k_n \in \AA$ with
$k_n \res H = (f_{n+1} - f_{n}) \res H$ and $\|k_n\| \le c 2^{-n}$.
Let $g = f_0 + \sum_{n=0}^\infty k_n$.
\end{proof}

If $f \in C(X)$, then $f$ is an \textit{idempotent} iff $f^2 = f$.
Note that $f$ is an idempotent iff 
$f = \cchi_H$ for some clopen $H \subseteq X$.

\begin{definition}
\label{def-equiv}
Assume that $\AA \lea C(X)$ and $\AA$ is closed.
\begin{itemizn}{"2B}
\item
$\BB_\AA$ is
the set of all clopen $H \subseteq X$ such that $\cchi_H \in \AA$.
\item
If $x,y\in X$, say that $x \sim_\AA y$ iff 
every $H \in \BB_\AA$ contains either both or neither of $x,y$.
We delete the subscript $\AA$ when it is clear from context.
\item
For $P \subseteq X$, $\widetilde P = \bigcap \{H \in \BB_\AA : P \subseteq H\}$.
\item
$P \subseteq H$ \emph{factors through} $\sim$ iff
$P$ is a union of $\sim$ -- equivalence classes.
\end{itemizn}
\end{definition}

Some simple observations:
$\BB_\AA$ is a boolean algebra, and each equivalence
class is closed and factors through $\sim$.
$\widetilde P$ is always closed; and
if $P$ is closed in $X$, 
then, by compactness, $P$ factors through $\sim$ iff
$P = \widetilde P$.
If $X$ is totally disconnected, then $\AA = C(X)$ iff 
$\BB_\AA$ contains all clopen subsets of $\AA$ iff
each equivalence class is a singleton.
When $P = \widetilde P$, Lemma \ref{lemma-restrict-factor} below
describes $\BB_{\AA\res P}$ in terms of $\BB_\AA$.
As in \cite{RUD1, HS}, we need the following lemma
for producing elements of $\BB_\AA$:

\begin{lemma}
\label{lemma-get-clopen}
Suppose that $\AA \lea C(X)$ and $\AA$ is closed, and suppose that
$h \in \AA$ and $b \in \RRR \setminus \Re(h(X))$.
Then $\{x \in X : \Re(h(x)) <  b\} \in \BB_\AA$.
\end{lemma}

As usual, $\Re(z)$ denotes the real part of the complex number $z$;
so $\Re: \CCC \to \RRR$.
Note that  $\{x \in X : \Re(h(x)) <  b\} $ is clearly clopen,
and the result is trivial unless $ \Re(h(X))$ contains elements
from both $(-\infty, b)$ and $(b, \infty)$.
In that case, the lemma is easily proved using Runge's Theorem.

It is easy to see that
$\AA\res H $ is closed in $C(H)$ for each $H \in \BB_\AA$.
More generally,

\begin{lemma}
\label{lemma-restrict-factor}
Assume that $\AA\lea C(X)$ and is closed in $C(X)$, and assume that
$P$ is closed in $X$ and factors through $\sim_\AA$.
Then
\begin{itemizz}
\item [1.] $\AA \res P$ is closed in $C(P)$.
\item [2.] $\BB_{\AA\res P} = \{H \cap P : H \in \BB_\AA\}$.
\item [3.] If $P$ is an equivalence class,
then $\BB_{\AA\res P} = \{\emptyset, P\}$.
\end{itemizz}
\end{lemma}
\begin{proof}
For (1):
We verify (2) of Lemma \ref{lemma-closed-equiv}, with $c = 2$.
Fix $f \in \AA$ with $f \res P$ not identically $0$.
Since $P = \widetilde P$,
we may choose $K \in \BB_\AA$ with $P \subseteq K$ and
$\|f\|_K \le 2 \|f\|_P$.  Now let $ f^* = f \cdot \cchi_K$.
Then $ f^* \in \AA$ and $\| f^*\| \le 2 \|f\|_P$, and
$f \res P =  f^*  \res P$.

For (2):  Fix $H \in \BB_{\AA\res P}$, and we produce an
$\widehat H \in \BB_\AA$ with $H = \widehat H \cap P$.
Fix $f \in \AA$ such that $f\res P = \cchi_H$.
Since $P = \widetilde P$,
we may choose $K \in \BB_\AA$ such that
$\min(|f(x)|, |f(x)-1|) < 1/4$ for all $x \in K$.
Let $g = f \cdot \cchi_K$.  Then $1/2 \notin \Re(g(X))$, so,
by Lemma \ref{lemma-get-clopen},
$\widehat H := \{x \in X : \Re(g(x)) >  1/2\} \in \BB_\AA$.

(3) is immediate from (2).
\end{proof}

Note that in (3), if the equivalence class $P$ is not
a singleton, then $P$ cannot be scattered by Theorem \ref{thm-rudin}.
$P$  may have isolated points (see Example \ref{ex-circ-dot} below),
but if $Q = \ker(P)$, then $\BB_{\AA\res Q} = \{\emptyset, Q\}$
(see Lemma \ref{lemma-restrict-kernel}).
Our proof of Theorem \ref{thm-lots} will focus on showing that
for appropriate perfect $Q$, this situation cannot happen --- i.e.,
$\AA\res Q$ must contain non-trivial idempotents.

We now define two important closed subsets of $X$ associated with $\AA$:

\begin{definition}
Assume that $\AA\lea C(X)$ and is closed in $C(X)$.
Let $H$ be a closed subset of $X$.  Then:
\begin{itemizn}{"2B}
\item
$H$ is \emph{essential} for $\AA$ iff
for all $f,g \in C(X)$ with $f \res H = g\res H$:
$f\in\AA$ iff $g \in \AA$.
\item
$E(\AA)$ denotes the \emph{essential set} for $\AA$;
this is the smallest closed set which is essential for $\AA$.
\item
$H$ is a \emph{boundary} for $\AA$ iff
$\|f\|_H = \|f\|$ for all $f\in \AA$.
\item
$\Sh(\AA)$ denotes the  \v Silov \emph{boundary};
this is the smallest closed set which is a boundary for $\AA$.
\end{itemizn}
\end{definition}
Bear \cite{BEAR1} defined the notion of ``essential''
and proved the existence of a smallest essential set.
The existence of a smallest boundary, $\Sh(\AA)$,
is a classical result of  \v Silov;
short proofs are given in \cite{BEAR2, GA, RUD4}.
If we fix $X$ and let $\AA$ vary, then note that as $\AA$ increases,
$E(\AA)$ decreases and $\Sh(\AA)$ increases.
$E(\AA) = \emptyset$ iff $\AA = C(X)$.
$\Sh(\AA)$ is never empty and $\Sh(C(X)) = X$.
Proposition \ref{prop-ess-sil} contains some additional properties 
of $E(\AA)$ and $\Sh(\AA)$.

\begin{example}
\label{ex-circ-dot}
\rm
Let $\AA_0 \lea C(\overline D)$ be the disc algebra.
Let $X = \TTT \cup H \cup F$, where $H,F \subset \CCC$ are compact scattered,
with $H \subset D$ and $F \cap \overline D = \emptyset$.
Let $\AA_1 = \AA \rest (\TTT \cup H)$.
Let $\AA = \{g\in C(X): g\res (\TTT \cup H) \in \AA_1\}$.
Then the $\sim_\AA$ -- equivalence classes
are $P := \TTT \cup H$ and the singletons in $F$.
$\ker(P) = \ker(X) = \TTT$,
$\AA \res \TTT$ is closed in $C(\TTT)$, and 
$\BB_{\AA \res \TTT} = \{\emptyset, \TTT \}$.
$E(\AA) = \TTT \cup H$ and $\Sh(\AA) = \TTT \cup F$.
$\Sh(\AA\res E(\AA)) = \TTT$; the fact that
$\Sh(\AA\res E(\AA))$ is perfect holds in general;
see Proposition \ref{prop-ess-sil}.
\end{example}

The next lemma shows that the conclusions (1)(2)
of Lemma \ref{lemma-restrict-factor}
hold for many closed subsets of $X$ which do
not necessarily factor through  $\sim_\AA$:

\begin{lemma}
\label{lemma-restrict-boundary}
Assume that $\AA\lea C(X)$ and is closed in $C(X)$, and assume that
$K$ is closed in $X$. 
If $\Sh(\AA) \subseteq K$ then {\rm(1)(2)} below hold.
If $E(\AA) \subseteq K$, then {\rm (1)} holds, and {\rm (2)} holds when
$X$ is totally disconnected.
\begin{itemizz}
\item [1.] $\AA \res K$ is closed in $C(K)$.
\item [2.] $\BB_{\AA\res K} = \{H \cap K : H \in \BB_\AA\}$.
\end{itemizz}
\end{lemma}
\begin{proof}
(1) is immediate from Lemma \ref{lemma-closed-equiv}.
For (2), when $\Sh(\AA) \subseteq K$,
fix $H \in \BB_{\AA\res K}$, and then fix
$f\in \AA$ such that $f \res K = \cchi_H$.
Then $f^2 - f$ is $0$ on $K$ and hence is $0$ everywhere,
so $f$ is an idempotent.  
If $\widehat H = f\iv\{1\}$, then $\widehat H \in \BB_\AA$
and $H = \widehat H \cap K$.
\end{proof}

\begin{lemma}
\label{lemma-sh-all}
Assume that $\AA\lea C(X)$ and is closed in $C(X)$, with $\AA\ne C(X)$.
Then $\AA\res E(\AA) \ne C(E(\AA))$ and  
$\AA\res \Sh(\AA) \ne C(\Sh(\AA))$.
\end{lemma}
\begin{proof}
The first statement is clear from the definition of $E(\AA)$.
The second is clear unless $\Sh(\AA) \ne X$, so fix
$p \in X \backslash \Sh(\AA)$.  
Define $\varphi: \AA\,\res\,\Sh(\AA) \to \CCC$ by
$\varphi(f \res \Sh(\AA)) = f(p)$.
Then $\varphi \in \MM(\AA\,\res\,\Sh(\AA))$ and
$\varphi$ differs from all point evaluations,
so $\AA\,\res\, \Sh(\AA) \ne C(\Sh(\AA))$.
\end{proof}

The next two lemmas relate $\Sh(\AA)$ and $E(\AA)$ to the equivalence classes:

\begin{lemma}
\label{lemma-sh-dense}
Let $\AA \lea C(X)$ be closed.
Then $\Sh(\AA)$ meets every $\sim_\AA$ -- equivalence class.
\end{lemma}

\begin{lemma}
\label{lemma-in-E}
Let $\AA \lea C(X)$ be closed.
Let $H$ be closed in $X$ such that $\{p\}$ is an
$\sim_\AA$ -- equivalence class for all $p \in X \backslash H$.
Then $E(\AA) \subseteq H$.
\end{lemma}

\section{Restricting to the Kernel}
\label{sec-ker}

We consider the restrictions of $\AA \lea C(X)$ to $K = \ker(X)$
and to $\widetilde K$ (see Definition \ref{def-equiv}).
In Example \ref{ex-circ-dot},
$\widetilde K$ properly contains $K$ (i.e., $K$ does not factor
through $\sim$), and $K$ is neither a boundary nor essential.
Nevertheless, the next two lemmas show that the conclusions (1)(2) of 
Lemmas \ref{lemma-restrict-factor} and \ref{lemma-restrict-boundary}
hold for $K$:

\begin{lemma}
\label{lemma-shkern}
Assume that $\AA\lea C(X)$ and is closed in $C(X)$.
Assume that $X$ is not scattered and let
$K = \ker(X)$.  Then
$\Sh(\AA \res \widetilde K) \subseteq K$.
\end{lemma}
\begin{proof}
If not, then $K$ is not a boundary for $\AA \res \widetilde K$,
so fix $f\in\AA$ such that $\|f\|_K < \|f\|_{\widetilde K}$.
Multiplying by a constant, we may assume that
$\|f\|_K \le 1$ and $f(p) = 1 + 2\varepsilon > 1$ for some $p \in \widetilde K$.
Let $W \subset X$ be clopen with $K \subseteq W$
and $f(W) \subseteq B(0; 1 + \varepsilon)$,
so that $\Re(f(W)) \subseteq (1 - \varepsilon,1 + \varepsilon )$.
Since $X\backslash W$
is scattered and compact, so is $\Re(f(X\backslash W))$, so we can fix
$b \in (1 + \varepsilon, 1 + 2\varepsilon)$ such that $b \notin \Re(f(X))$.
Let
$H = \{x \in X : \Re(f(x)) <  b\}$.
Then $K \subseteq H$ and $p \notin H$, and, 
by Lemma \ref{lemma-get-clopen}, $H \in \BB_\AA$, contradicting
$p \in \widetilde K$.
\end{proof}

\begin{lemma}
\label{lemma-restrict-kernel}
Assume that $\AA\lea C(X)$ and is closed in $C(X)$, and let
$K = \ker(X)$.  Then
\begin{itemizz}
\item [1.] $\AA \res K$ is closed in $C(K)$.
\item [2.] $\BB_{\AA\res K} = \{H \cap K : H \in \BB_\AA\}$.
\end{itemizz}
\end{lemma}
\begin{proof}
By Lemmas \ref{lemma-shkern}, 
\ref{lemma-restrict-factor},  and
\ref{lemma-restrict-boundary}.
\end{proof}

\begin{definition}
The compact space $X$
has the \emph{NTIP}
iff every closed $\AA \lea C(X)$ has a non-trivial idempotent
$(\!$i.e., of the form $\cchi_H$, with
$\emptyset \subsetneqq H \subsetneqq X)$.
\end{definition}

So, the NTIP is trivially false of connected spaces.
If $X$ is not connected, then the CSWP implies the NTIP\@.
There are totally disconnected $X$ which have the NTIP
but not the CSWP; see Example \ref{ex-ntip}.
Lemma \ref{lemma-restrict-kernel} implies immediately:

\begin{lemma}
\label{lemma-ntip-ker}
For any compact non-scattered $X$, if
$\ker(X)$ has the NTIP then $X$ has the NTIP.
\end{lemma}

The converse is false even for totally disconnected $X$;
see Example \ref{ex-ntip}.

\begin{lemma}
\label{lemma-restrict-to-perfect}
Assume that $X$ is compact 
and every perfect subset of $X$ has the NTIP\@.  Then $X$ has the CSWP.
\end{lemma}
\begin{proof}
Let $\AA \lea C(X)$ be closed but not all of $C(X)$.
Then there must be a
$\sim_\AA$ -- equivalence
class $P$ which is not a singleton.  By Lemma
\ref{lemma-restrict-factor}, $P$ does not have the NTIP, so $P$ is not
scattered.  If $Q = \ker(P)$, then $Q$ is perfect and does
not have the NTIP by Lemma \ref{lemma-ntip-ker}.
\end{proof}

\begin{lemma}
\label{lemma-dense-ker}
Let $\AA \lea C(X)$, assume that $X$ is not scattered,
and let $K = \ker(X)$.
If $\AA \res K$ is dense in $C(K)$,
then $\AA$ is dense in $C(X)$.
\end{lemma}
\begin{proof}
WLOG, $\AA$ is closed in $C(X)$, so
$\AA \res K = C(K)$ by Lemma \ref{lemma-restrict-kernel}.
By Lemma \ref{lemma-shkern}, we have
$\Sh(\AA \res \widetilde K) \subseteq K$, so that
$\AA \res \widetilde K = C(\widetilde K)$ by Lemma \ref{lemma-sh-all},
which implies that $\widetilde K = K$.
It follows now by Rudin's Theorem \ref{thm-rudin}
that each point in $X \backslash K$ is a $\sim$ equivalence class,
so $E(\AA) \subseteq K$ by Lemma \ref{lemma-in-E}.  This, plus
$\AA \res K = C(K)$, implies that $\AA = C(X)$.
\end{proof}

\begin{corollary}
\label{cor-CSWP-iff}
$X$ has the CSWP iff $\ker(X)$ has the CSWP.
\end{corollary}

\section{Compact Ordered Spaces}
\label{sec-lots}
A \textit{linearly ordered topological space}, or LOTS, 
is a topological space $X$ with a total order $<$ such that
the topology on $X$ is the order topology generated by $<$.
We begin with some notation and well-known facts about
such spaces.

\begin{definition}
\label{def-neighbors}
In any set $X$ totally ordered by $<$:
\begin{itemizn}{"2B}
\item
$0_X$ and $1_X$ {\rm(}or $0,1$, when $X$ is clear from context\/{\rm)}
denote,
respectively, the first and last elements of $X$ {\rm(}when these
elements exist\/{\rm)}.
\item
$x^- \ltn x^+$ holds iff $x^- < x^+$ and there are no
points of $X$ strictly between $x^-$ and $x^+$.
In that case, we call $\{x^-, x^+\}$ a \emph{neighboring pair},
and say that $x^-$ and  $x^+$ are \emph{neighbors}.
\end{itemizn}
\end{definition}

If $X$ is a compact LOTS with no isolated points, then
the points $0$ and $1$ exist and have no neighbors; and
each $x \in X \backslash \{0,1\}$ has either no neighbors
(when it is a limit point from the left and the right),
or exactly one neighbor (when it is a limit point from 
one side and not the other).

The standard unit interval in the reals is connected,
so there are no neighboring pairs, but one can form examples with
neighboring pairs by doubling points in the interval:

\begin{definition}
If $E \subseteq (0,1) \subseteq \RRR$,
then $\DDD(E) = [0,1] \times \{0\} \cup E \times \{1\}$, ordered
lexicographically.
\end{definition}

Equivalently, we form $\DDD(E)$ by doubling every $x$ in $E$
(with $x^-$ represented by $(x,0)$ and $x^+$ represented by $(x,1)$),
and not doubling the points of $[0,1] \setminus E$.  We never
double $0$ or $1$, since that would create an isolated point.

\begin{proposition}
\label{prop-cantor}
For every $E \subseteq (0,1) \subseteq \RRR$,
the space $\DDD(E)$ is a compact separable LOTS with no
isolated points.  $\DDD(E)$ contains a copy of the Cantor
set iff in the reals, $(0,1)\backslash E$ contains a copy of the Cantor set.
\end{proposition}

The \textit{double arrow space} is $\DDD( (0,1) )$.
The proof of Theorem \ref{thm-lots} will involve only
spaces of the form $\DDD(E)$, by the following four lemmas.
The first is from Lutzer and Bennett \cite{LB}.

\begin{lemma}
\label{lemma-hs}
If $X$ is a separable LOTS and $H \subseteq X$, then $H$ is separable
and Lindel\"of in its relative topology inherited from $X$.
\end{lemma}

The relative topology on $H$ is not in general the LOTS topology
induced by the order restricted to $H$, but

\begin{lemma}
\label{lemma-rel}
If $X$ is a compact LOTS and $H$ is a closed
subset of $H$, then the relative topology and the order topology
agree on $H$.
\end{lemma}

\begin{lemma}
\label{lemma-disc}
If $X$ is a compact separable LOTS 
and $X$ does not contain a copy of the Cantor set, then 
$X$ is totally disconnected.
\end{lemma}
\begin{proof}
If $a < b$ and  $[a,b] \subseteq X$ is connected, then
$[a,b]$ is isomorphic to a closed
interval in $\RRR$, and hence contains a Cantor subset.
\end{proof}

\begin{proposition}
\label{prop-standard}
If $X$ is a compact separable LOTS with no isolated points,
then $X$ is isomorphic to $\DDD(E)$ for some
$E \subseteq (0,1) \subseteq \RRR$.
\end{proposition}

To prove Theorem \ref{thm-lots}, we need to show
that every compact separable LOTS $Y$ which does not contain
a copy of the Cantor set has the CSWP\@.  By
Lemma \ref{lemma-restrict-to-perfect}, it is enough to show
that if $X$ is any perfect subset of $Y$, then $X$
has the NTIP\@.  By Lemmas \ref{lemma-hs} and \ref{lemma-rel},
$X$ is also a compact separable LOTS and does not contain a Cantor subset.
By Lemma \ref{lemma-disc}, $X$ and $Y$ are totally disconnected.
Propositions \ref{prop-cantor} and
\ref{prop-standard} are not used in the proof of Theorem \ref{thm-lots},
but they characterize the spaces to which Theorem \ref{thm-lots} applies.

Lemma \ref{lemma-ntip} will
establish the NTIP for $X$.  It is proved by analyzing step functions.
In general, $\sigma \in C(X)$ is \textit{step function} 
iff $\sigma(X)$ is finite.  If $X$ is compact and totally disconnected,
then the step functions are dense in $C(X)$.
If $X$ is also a LOTS, then there is a simple description of
such step functions:

\begin{definition}
\label{def-step}
If $X$ is a compact totally disconnected LOTS with no isolated
points and $S \subseteq X$, then:
\begin{itemizn}{"2B}
\item
A \emph{step function with endpoints in $S$} is a 
function of the form $\sum_{i = 0}^n z_i \cchi_{I_i}$,
where $n \ge 0$, each $z_i \in \CCC$, and each $I_i = [a_i^+, a_{i+1}^-]$
is a clopen interval in $X$, with
$0 = a_0^+ < a_1^-  \ltn  a_1^+   <   \cdots <
       a_{n}^-  \ltn      a_{n}^+ <   a_{n+1}^- = 1$, where all
$a_0^+ , a_1^- , a_1^+ ,\cdots, a_{n}^- , a_{n}^+ , a_{n+1}^- \in S$.
\item
$S$ is \emph{nice} iff $S$ is of the form
$\{0,1\} \cup \bigcup_\alpha \{x_\alpha^-,x_\alpha^+\}$, where each
$x_\alpha^- \ltn x_\alpha^+$.
\item
$\SF(S)$ is the set of all step functions with endpoints in $S$.
\end{itemizn}
\end{definition}

When studying step functions with endpoints in $S$,
we consider only nice $S$ because only $0,1$, and neighboring
pairs get used as endpoints.
Note that $\SF(S)$ is always a sub-algebra of $C(X)$.
$\SF(\{0,1\})$ is the set of constant functions (the $n=0$
case of Definition \ref{def-step}).
$\SF(S)$ is dense in $C(X)$
in the case that $S$ contains all neighboring pairs.
However, if $S$ omits some neighboring pair $\{x^-,x^+\}$, then
$\cchi_{[0,x^-]}$ is a step function which has distance
$1/2$ from $\SF(S)$.

\begin{definition}
For $X$ a compact LOTS, $f \in C(X)$ and $\varepsilon > 0$, 
\[
\JMP_\varepsilon(f) = \{0,1\} \cup
\bigcup\{ \{x^-, x^+\} : x^- \ltn x^+ \ \&\  |f(x^-)- f(x^+)|
\ge \varepsilon \}\ \ .
\]
$\JMP(f) = \bigcup\{\JMP_\varepsilon : \varepsilon > 0\}$.
\end{definition}

So, $\JMP(f)$ contains all neighboring pairs at which $f$ \textit{jumps}.
By compactness of $X$ and continuity of $f$,

\begin{lemma}
\label{lemma-J-finite}
Each $\JMP\varepsilon(f)$ is finite, and $\JMP(f)$ is countable.
\end{lemma}

Now, suppose that the $X$ of Definition \ref{def-step} is separable
but not second countable.  Then we may choose a nice $S$ with
$S$ countable and dense in $X$, but there are uncountably many
neighboring pairs, so 
$\SF(S)$ will not be dense in $C(X)$. 
Specifically, say $f\in C(X)$, $\varepsilon > 0$,
and $x^\pm \in \JMP_{2\varepsilon}(f) \setminus S$.
Then there is a step function
$\sigma$ with $\|f - \sigma\| < \varepsilon$,
but the endpoints $a_i^\pm$ will not all lie in $S$,
since $x^\pm$ is among the $a_i^\pm$.
However, we can choose $b_i^\pm \in S$ with
$ a_i^+ < b_i^-  \ltn  b_i^+  <   a_{i+1}^- $, and then ``describe''
$\sigma$ just using points in $S$ as follows:

\begin{definition}
\label{def-desc}
If $X$ is a compact totally disconnected LOTS with no isolated
points and $S \subseteq X$, then
a \emph{step function descriptor from} $S$ is
a finite array of the form:
\[
\Delta = 
\left(
\begin{array}{cccc}
     z_0     &       z_1    & \cdots\cdots & z_{n}    \\
b_0^-, b_0^+ & b_1^-, b_1^+ &\cdots\cdots  &  b_{n}^-, b_{n}^+
\end{array}
\right)
\ \ ,
\]
where
$0 < b_0^- \ltn b_0^+ < b_1^-\ltn b_1^+ < \cdots < b_{n}^- \ltn b_{n}^+ < 1$
in $X$, with each  $b_i^\pm \in S$ and each $z_i \in \CCC$.
Then $\STEP(\Delta)$ is the set of all step functions $\sigma$ which
\emph{meet this description} in the sense that
$\sigma = \sum_{i = 0}^n z_i \cchi_{I_i}$ for some
$I_i = [a_i^+, a_{i+1}^-]$, where
\[
0 = a_0^+ < b_0^- \ltn  b_0^+   <  a_1^-  \ltn  a_1^+   <  b_1^+
 \ltn  b_1^+  <  \cdots <
   a_{n}^-  \ltn  a_{n}^+ <  b_{n}^+  \ltn  b_{n}^+  < a_{n+1}^- = 1
\ \ .
\]
\end{definition}

Note that the endpoints of $\sigma$ are not required to be in $S$.

\begin{lemma}
\label{lemma-desc2}
Suppose $\Delta$ is as in Definition \ref{def-desc}, 
$\sigma, \tau \in \STEP(\Delta)$, $\varepsilon > 0$,
and $f, g \in C(X)$ with
$\|f - \sigma\| < \varepsilon$ and
$\|g - \tau \| < \varepsilon$.  Let $h = f - g$.  Then
\[
h(X) \subseteq B(0; 2\varepsilon) \cup
\bigcup_{i < n} B(w_i; 2\varepsilon) \ \ ,
\]
where each $w_i$ is either $z_i - z_{i+1}$ or $z_{i+1}  -  z_i$.
\end{lemma}

With appropriate values of the $z_i$, and $\varepsilon$ small enough,
this will force $\Re(h(X))$ to be disconnected, so that we can
apply Lemma \ref{lemma-get-clopen} to get a non-trivial idempotent in $\AA$.
We shall prove that $\Re(h(X))$  is disconnected by using:

\begin{lemma}
\label{lemma-proj-disc}
Assume that
$r \ge 1$,
$F \subseteq \CCC$, 
$F \cap B(0;1) \ne \emptyset$,
$F \cap B(w_0;1/(3r)) \ne \emptyset$,
and $F \subseteq B(0;1) \cup \bigcup_{k < r}B(w_k; 1/(3r))$
for some $w_0, \ldots, w_{r-1} \in \CCC$ with $w_0 = \pm 2$.
Then $\Re(F)$ is not connected.
\end{lemma}

We shall get $h = f - g \in \AA$ by applying the following lemma:

\begin{lemma}
\label{lemma-lst}
Assume that $X$ is a compact separable totally disconnected LOTS
with no isolated points and $\AA \lea C(X)$.
Then there is a countable dense nice 
$S \subseteq X$ such that for all step function descriptors $\Delta$
from $S$ and all $\varepsilon > 0$:
\begin{itemizn}{"63}
\item
If there is a $\sigma \in \STEP(\Delta)$ and an $f \in \AA$
with $\|f - \sigma\| < \varepsilon$,  then there is a step function
$\tau \in \STEP(\Delta)$ \emph{with endpoints in} $S$ and
some $g \in \AA$
with $\|g - \tau\| < \varepsilon$.
\end{itemizn}
\end{lemma}
\begin{proof}
Call a step function descriptor $\Delta$ \textit{rational} iff
the $z_i$ from Definition \ref{def-desc} all have rational real and 
imaginary parts.  Observe that in proving the lemma, it is sufficient
to obtain \ding{"63} for all rational $\varepsilon$ and
all rational $\Delta$.
Let $S_0 \subseteq X$ be nice and countable and dense.
Then, obtain an increasing chain
$S_0 \subseteq S_1 \subseteq S_2 \cdots$ as follows:

Given $S_n$, and given a rational step function descriptor $\Delta$
from $S_n$ and a rational $\varepsilon > 0$:
If there exists a $\sigma \in \STEP(\Delta)$ and an $f \in \AA$
with $\|f - \sigma\| < \varepsilon$,  then let
$\sigma^n_{\Delta, \varepsilon} \in \STEP(\Delta)$ and
$f^n_{\Delta, \varepsilon} \in \AA$ be some such $\sigma,f$.
Let $S_{n+1}$ be $S_n$ together with all the endpoints
of all the $\sigma^n_{\Delta, \varepsilon}$.

Let $S = \bigcup_{n\in\NNN} S_n$.
\end{proof}

Putting these together, we get:

\begin{lemma}
\label{lemma-ntip}
Assume that $X$ is a compact totally disconnected LOTS which is
separable but not second countable.
Then $X$ has the NTIP.
\end{lemma}
\begin{proof}
Note that by Lemma \ref{lemma-hs},
$\ker(X)$ is also separable but not second countable,
and it is sufficient to prove that $\ker(X)$ has the NTIP by
Lemma \ref{lemma-ntip-ker}.  Thus, we may assume that $X$ is perfect.

Let $\AA \lea C(X)$ be closed.

First, fix $S  \subseteq X$ as in Lemma \ref{lemma-lst}.
Then, since $X$ is not second countable, we can fix
a neighboring pair $q^- \ltn q^+$ with $q^\pm \notin S$.
Then fix an $f\in \AA$ with $f(q^+) = 1$ and $f(q^-) = -1$.

Let $R = \JMP_{1/3}(f)$, which is finite by Lemma \ref{lemma-J-finite},
and let $|R| = 2r + 2$; that is, $R$ contains
$0$ and $1$ plus $r$ neighboring pairs.  Then $r \ge 1$, since $q^\pm \in R$.

Fix a step function $\sigma$ with $\|f - \sigma\| < 1/(6 r)$.
Say $\sigma = \sum_{i = 0}^n z_i \cchi_{I_i}$,
where  each $z_i \in \CCC$, each $I_i = [a_i^+, a_{i+1}^-]$, and
$0 = a_0^+ < a_1^-  \ltn  a_1^+   <   \cdots <
       a_{n}^-  \ltn      a_{n}^+ <   a_{n+1}^- = 1$.
Then $n \ge r$, since the $a_i^\pm$ include in particular
the neighboring pairs from $R$.

Choose the $b_i^\pm \in S$ and $\Delta$ as in Definition
\ref{def-desc} so that $\sigma \in \STEP(\Delta)$, and then
apply  \ding{"63} from  Lemma \ref{lemma-lst} and choose
$\tau \in \STEP(\Delta)$ with endpoints in $S$ and $g \in \AA$
with $\|g - \tau\| < 1/(6r)$.
We now have
\[
{
\setlength{\arraycolsep}{2pt}
\begin{array}{cccccccccccccccccccc}
0 &=& a_0^+ &<& b_0^- \ltn  b_0^+   &<&  a_1^-  \ltn  a_1^+   &<&  b_1^+
 \ltn  b_1^+  &<&  \cdots &<&
   a_{n}^-  \ltn  a_{n}^+ &<&  b_{n}^+  \ltn  b_{n}^+  &<& a_{n+1}^- &=& 1&\\
0 &=& c_0^+ &<& b_0^- \ltn  b_0^+   &<&  c_1^-  \ltn  c_1^+   &<&  b_1^+
 \ltn  b_1^+  &<&  \cdots &<&
   c_{n}^-  \ltn  c_{n}^+ &<&  b_{n}^+  \ltn  b_{n}^+  &<& c_{n+1}^- &=& 1&\ \ ,
\end{array}
}
\]
where the $c_i^\pm \in S$ and  $\tau = \sum_{i = 0}^n z_i \cchi_{J_i}$,
with each $J_i = [c_i^+, c_{i+1}^-]$.  Let $h = f - g$.
By Lemma \ref{lemma-desc2}, we have
\[
F := h(X) \subseteq B(0;  1/(3r)  ) \cup
\bigcup_{i < n} B(w_i;  1/(3r) ) \ \ ,
\]
where each $w_i$ is either $z_i - z_{i+1}$ or $z_{i+1}  -  z_i$.

Next, note that $|w_i| \le 2/3$ for all but at most $r$ indices $i$:
This holds because $\|f - \sigma\| < 1/6$, so if $|w_i| > 2/3$, then
$|f(a_i^-) - f(a_i^+)| > 1/3$, so that $a_i^\pm$ is one of the
$r$ neighboring pairs in $R$.

$F \cap B(0; 1) \ne \emptyset$ because $h(0) \in B(0; 1)$.
Also, one of the $a_i^\pm$ is the pair $q^\pm$, and this 
$w_i = \pm 2$.  Since this $a_i^\pm = q^\pm \ne S$,
while $c_i^\pm \in S$, $\; F \cap B(w_i; 1/(3r)) \ne \emptyset$.

Furthermore, whenever $|w_i| \le 2/3$, we have
$B(w_i;  1/(3r) )  \subseteq B(0;1)$, and this holds for
all but at most $r$ of the indices $i$.  It follows now by
Lemma \ref{lemma-proj-disc} that $\Re(h(X))$  is disconnected,
so that by Lemma \ref{lemma-get-clopen},
$\AA$ contains a non-trivial idempotent.
\end{proof}

This lets us prove Theorem \ref{thm-lots}
in the separable case:

\begin{lemma}
\label{lemma-lots-separable}
Let $X$ be totally ordered by $<$, and assume that $X$
is compact and separable in its order topology.  Then $X$ has the CSWP
iff $X$ does not contain a copy of the Cantor set.
\end{lemma}
\begin{proof}
To prove that $X$ has the CSWP, it is sufficient, 
by Lemma \ref{lemma-restrict-to-perfect}, 
to prove that every perfect $P \subseteq X$ has the NTIP, 
and this holds by Lemma \ref{lemma-ntip}.
\end{proof}

The reader familiar with elementary submodel techniques
(see Dow \cite{DOW}) will note that one could simplify somewhat
the proof of Lemma \ref{lemma-ntip},
at the expense of introducing some new concepts into the proof.
Specifically, 
if $\theta$ is a suitably large regular cardinal
and $\AA,X \in M \prec H(\theta)$, with $M$ countable,
then we can simply define $S$ in Lemma \ref{lemma-ntip}
to be $M \cap X$, and omit Lemma \ref{lemma-lst} entirely.

\section{Measures and Supports}
\label{sec-meas}

\begin{definition}
If $\mu$ is a regular complex Borel measure on the compact space $X$,
then $|\mu|$ denotes its total variation, and
$\supt(\mu) = \supt(|\mu|)$ denotes its (closed) support;
that is,
$\supt(\mu) = X \setminus \bigcup \{U\subseteq X :\; $U$ \mbox{ is open }
\ \ \& \ \ |\mu|(U) = 0\}$.
\end{definition}

\begin{lemma}
\label{lemma-supt-cswp}
Assume that $X$ is compact and that $\supt(\mu)$ has the CSWP
for all regular Borel measure $\mu$.
Then $X$ has the CSWP.
\end{lemma}
\begin{proof}
Suppose we have $\AA \lea C(X)$ with $\AA$ closed in $C(X)$ and
$\AA \ne C(X)$.
Just viewing $\AA$ and $C(X)$ as Banach spaces, fix $\varphi \in C(X)^*$
such that $\varphi(f) = 0$ for all $f \in \AA$, but such that
$\varphi(g) \ne 0$ for some $g \in C(X)$.
Let $\mu$ be a regular complex Borel measure  such
that $\varphi(f) = \int f \, d \mu$ for all $f \in C(X)$.
Let $H = \supt(\mu)$.

Then $\int_H h \, d \mu = 0$ for all $h \in \AA \res H$,
but $\int_H (g\res H) \, d \mu = \varphi(g) \ne 0$, so $\AA \res H$
is not dense in $C(H)$.  Thus, $H$ does not have the CSWP.
\end{proof}

In the case of a LOTS, it is easy to see
(by a minor modification of the proof of Theorem 1.0 in \cite{DK}):

\begin{lemma}
\label{lemma-lots-supt}
Assume that $X$ is a compact LOTS and $\mu$ is a regular Borel measure
on $X$.  Then $\supt(\mu)$ is separable.
\end{lemma}

\begin{proof}[\bf Proof of Theorem \ref{thm-lots}]
Immediate from Lemmas
\ref{lemma-lots-separable},
\ref{lemma-supt-cswp}, and
\ref{lemma-lots-supt}.
\end{proof}

\section{Remarks and Examples}
\label{sec-ex}

Much is known about the maximal ideal space,
$\MM(\AA)$, for some standard commutative unital Banach
algebras $\AA$.  For example, 
it is well-known that $\MM(L^\infty(\TTT))$ is $\st(\BB)$,
where $\BB$ is the measure algebra of $\TTT$ (the measurable sets
modulo the null sets) and $\st(\BB)$ denotes its Stone space.

Hoffman and Singer \cite{HS} (see also \cite{GA,HO})
discuss $\MM(H^\infty)$,
where $H^\infty$ consists of those 
elements of $L^\infty(\TTT)$ whose negative Fourier coefficients vanish.
One may view $H^\infty$ as a closed subalgebra of 
$C(\MM(H^\infty))$ via the Gel'fand Transform (which is an isometry
in this case).  \cite{HS} shows, among other things, that
with this identification, $\Sh(H^\infty)$ is homeomorphic
to $\MM(L^\infty(\TTT)) = \st(\BB)$, and 
hence $H^\infty \lea C(\st(\BB))$, establishing that
$\st(\BB)$ fails to have the CSWP\@.
Since each of $\st(\BB)$ and $\beta\NNN$ contains the other,
it follows, as mentioned in the Introduction, that 
$\beta\NNN$ (and hence every compact space containing $\beta\NNN$)
fails to have the CSWP\@.
$\MM(H^\infty)$ also contains an open copy of $D$, and
by \cite{HS}, $\MM(H^\infty)$ properly contains $D \cup \Sh(H^\infty)$.
By Carleson's Corona Theorem (see \cite{GA}), $D$ is 
dense in $\MM(H^\infty)$; hence, $E(H^\infty) = \MM(H^\infty)$,
since $E(H^\infty) \supseteq D$ follows from that fact that each
function in $H^\infty$ is holomorphic on $D$.

We also have
$C(\TTT) \hookrightarrow L^\infty(\TTT)$.  That is,
with the natural identification of $f \in C(\TTT)$ with
its equivalence class $[f] \in L^\infty(\TTT)$,
we identify $C(\TTT)$ with closed subalgebra of $L^\infty(\TTT)$.

Between $C(\TTT)$ and $L^\infty(\TTT)$, there is the algebra
of functions which have everywhere a left and right limit.
Note that the usual notion of left and right limits make sense on
the circle.  If $g : \TTT \to \CCC$ and $z = e^{i\theta}\in\TTT$,
then $\lim_{w \to z+} g(w) $ means 
$\lim_{\varphi \to \theta+} g(e^{i\varphi}) $, and
$\lim_{w \to z-} g(w) $ means 
$\lim_{\varphi \to \theta-} g(e^{i\varphi}) $.

\begin{definition}
$J(\TTT)$ is the set of all $g : \TTT \to \CCC$ such that for all $z$:
$\lim_{w \to z+} g(w)$ and $\lim_{w \to z-} g(w)$ exist and
$g(z) = ( \lim_{w \to z+} g(w)  + \lim_{w \to z-} g(w) )/2$.
\end{definition}

We require $g(z)$ to be the average of its left and right limits
so that 
if $g,h \in J(\TTT)$ and $g \ne h$, then 
$\{z: g(z) \ne h(z)\}$ contains an interval, so
$g,h$ define different elements of $L^\infty$.
We now have $C(\TTT) \subset J(\TTT) \hookrightarrow L^\infty(\TTT)$.
Note that if  $f \in J(\TTT)$, then the left and right limits
of $f$ must be equal except at a countable set of points.
$ J(\TTT)$ is the closure (in the supremum norm) of the algebra of piecewise
continuous functions; these are the $f \in J(\TTT)$ whose
left and right limits are equal except at a finite set of points.
We can identify $\MM(J(\TTT))$ as the double arrow space by:

\begin{lemma}
\label{lemma-double-arrow}
$J(\TTT)$ is isometric with $C( \DDD( (0,1) ) )$,
so that $\MM(J(\TTT))$ is homeomorphic to $\DDD( (0,1) )$.
\end{lemma}
\begin{proof}
Define the map $\Psi : C( \DDD( (0,1) ) ) \to J(\TTT)$ so
that for $f\in C( \DDD( (0,1) ) )$ and $z \in \TTT$, we
compute $(\Psi(f))(z)$ as follows:
Let $z = e^{2\pi i x}$, with $x \in [0,1]$.
If $x \in (0,1)$, we have $x^\pm \in \DDD( (0,1) ) )$,
and we let $(\Psi(f))(z) = (f(x^+) + f(x^-))/2$.
If $x \in \{0,1\}$, so that $z = 1$, we let
$(\Psi(f))(1) = (f(0) + f(1))/2$.
\end{proof}

Thus, by Theorem \ref{thm-lots},

\begin{corollary}
\label{cor-jump}
Suppose that $\AA$ is a subalgebra of $J(\TTT)$ which contains the constant
functions and separates the points
of $\TTT$, and assume that for each $z\in\TTT$, some function in  $\AA$ has
a discontinuity at $z$.  Then
$\AA$ is dense in $J(\TTT)$.
\end{corollary}

In particular,
$J(\TTT) \cap H^\infty =  C(\TTT) \cap H^\infty$,
since otherwise,
by rotational symmetry,
$J(\TTT) \cap H^\infty$ would contradict Corollary \ref{cor-jump}.
However, this special case of \ref{cor-jump}
can be seen directly by integrating
with the Poisson kernel; in fact, by
Zalcman \cite{ZA}, if $f \in H^\infty$ is essentially 
discontinuous at a point $z \in \TTT$,
then its essential range at $z$ is uncountable, 
so in particular it cannot have a simple jump discontinuity.

\bigskip

The following lemma and corollary describe a class of
spaces with the NTIP:

\begin{lemma}
\label{lemma-ntip-big}
Suppose that the compact $X$ is not second countable
but $\ker(X)$ is second countable.  Then $X$ has the NTIP.
\end{lemma}
\begin{proof}
If $\ker(X)$ is empty, then $X$ is scattered, so $X$ has the CSWP
by Rudin's Theorem \ref{thm-rudin}.2.  So, assume that $X$ is not
scattered.

Let $\AA \lea C(X)$ be closed.
Then $\AA\res\ker(X)$ is closed by Lemma \ref{lemma-restrict-kernel}.
Furthermore, $\AA$ is not separable (since it separates points in $X$,
which is not second countable), while $\AA\res\ker(X)$ is separable,
so the restriction map from $\AA$ to $\AA\res\ker(X)$ is not 1-1.
Fix a non-zero $f\in \AA$ such that its restriction to $\ker(X)$ is zero.
We may assume that $\Re(f(X))$ contains more than just $0$.
But also, $\Re(f(X))$ is scattered (since it is a continuous image
of the compact scattered $X/\ker(X)$), and hence
$\Re(f(X))$ is disconnected, so that $\AA$ contains a non-trivial
idempotent by Lemma \ref{lemma-get-clopen}.
\end{proof}

\begin{corollary}
\label{cor-ntip-equiv}
Suppose that $X$ is totally disconnected, with
$\ker(X)$ second countable and non-empty.  Then $X$ has the NTIP
iff $X$ is not second countable.
\end{corollary}
\begin{proof}
One direction is by Lemma \ref{lemma-ntip-big}.
For the other direction:
Note that $\ker(X)$ is homeomorphic to the Cantor set.
Assume that $X$ is second countable, and follow \cite{RUD1}:
We may assume that $X \subset \CCC \subset S^2$,
the Riemann sphere, and that all non-empty (relatively) open
subsets of $\ker(X)$ have positive Lebesgue measure.
Let $\AA \lea C(S^2)$ consist of the functions holomorphic
on $S^2 \backslash \ker(X)$.  Then $f(\ker(X)) = f(X) = f(S^2)$ for
all $f \in \AA$, so all $f(X)$ are connected.
\end{proof}

\begin{example}
\label{ex-ntip}
There are compact totally disconnected $X,Y$ such that:
\begin{itemize}
\item[1.] $X$ has the NTIP but $\ker(X)$ does not have the NTIP.
\item[2.] $Y$ has the NTIP but not the CSWP.
\end{itemize}
\end{example}
\begin{proof}
For (1), use Corollary \ref{cor-ntip-equiv} with any $X$ such that
$X$ is not second countable but $\ker(X)$ is the Canter set;
for example, $X$ can be the Aleksandrov duplicate of the Cantor set.

For (2), $Y$ can be the same $X$.  Or, $Y$ can be the disjoint
sum of the Cantor set and the double arrow space.  This is a compact
LOTS and has the NTIP by Lemma \ref{lemma-ntip}.
\end{proof}

\medskip

Some additional properties of the \v Silov boundary and the
essential set are given by:

\begin{proposition}
\label{prop-ess-sil}
Assume that $\AA \lea C(X)$ and $\AA$ is closed in C(X):
\begin{itemizz}
\item[1.] $E(\AA) \cup \Sh(A) = X$.
\item[2.] $E(\AA) \cap \Sh(A) \ne \emptyset$ unless
$E(\AA) = \emptyset$ (equivalently, unless $\AA = C(X)$).
\item[3.] If $p \in \Sh(\AA)$ and is isolated in $\Sh(\AA)$, then $p$ is
isolated in $X$ and $\cchi_{\{p\}} \in \AA$,
so $p \ne E(\AA)$.
\item[4.] If $E(\AA) \ne \emptyset$, then $\Sh(\AA\res E(\AA))$ is perfect.
\end{itemizz}
\end{proposition}
\begin{proof}
For (1), if $p \notin E(\AA) \cup \Sh(A)$, fix $f\in C(X)$ with
$f(E(\AA) \cup \Sh(A)) = \{0\}$ and $f(p) = 1$.
Then $f\in\AA$ because $E(\AA)$ is essential, 
but $f$ contradicts the fact that $\Sh(\AA)$ is a boundary.

For (2), if  $E(\AA) \cap \Sh(A) = \emptyset$ and $E(\AA) \ne \emptyset$,
fix $f\in C(X)$ with
$f(E(\AA))= \{1\}$ and $f(\Sh(\AA))= \{0\}$.
Then $f$ yields a contradiction as in (1).

For (3), assume that $p$ is an isolated point of $\Sh(\AA)$.
Then $H := \Sh(\AA) \backslash \{p\}$ is closed and is not a boundary, so fix
$f\in\AA$ such that $\|f\| = 1$ but $\|f\|_H < 1$.
Then $|f(p)| = 1$, since $\Sh(\AA)$ is a boundary; multiplying
by a constant, we may assume that $f(p) = 1$.
Then $f^n \to \cchi_{\{p\}}$ on $\Sh(\AA)$ as $n \to \infty$.
Since $\AA \res \Sh(\AA)$ is closed in $C(\Sh(\AA))$ (see Lemma
\ref{lemma-restrict-boundary}), there is a $g \in \AA$
such that $g\res \Sh(\AA) = \cchi_{\{p\}}$.
Then $g^2 -g$ is $0$ on $\Sh(\AA)$, and hence everywhere,
so $g = \cchi_K$ for some clopen $K$, with $p \in K$ and
$H \cap K = \emptyset$.
If $K = \{p\}$, then $\cchi_{\{p\}} \in \AA$.
If $K \ne \{p\}$, fix $q \in K \backslash \{p\}$, and then fix 
$h\in \AA$ with $h(q) = 1$ and $h(p) = 0$;
then $h \cdot \cchi_K$ contradicts the fact
that $H \cup \{p\}$ is a boundary.

For (4), note that $E(\AA\res E(\AA)) = E(\AA)$, so that
applying (3) to $\AA\res E(\AA)$ yields that no
$p \in \Sh(\AA\res E(\AA))$ is isolated in
$\Sh(\AA\res E(\AA))$.
\end{proof}

(4) immediately implies Rudin's Theorem \ref{thm-rudin}.2,
and (4) yields another proof of Lemma \ref{lemma-dense-ker} (using
Lemmas \ref{lemma-restrict-boundary} and \ref{lemma-sh-all}).

\end{document}